\documentclass[10pt]{article}
\usepackage{amssymb,amsmath,graphicx,acronym,amsfonts}

\newtheorem{corollary}{Corollary}[section]
\newtheorem{theorem}{Theorem}[section]
\newtheorem{lemma}{Lemma}[section]
\newtheorem{definition}{Definition}[section]
\newtheorem{proposition}{Proposition}[section]
\newtheorem{example}{Example}[section]
\newtheorem{assum}{Assumption}[section]
\newtheorem{algo}{Algorithm}[section]
\newtheorem{Remark}{Remark}[section]

\def\bc{\begin{corl}}
\def\bc{\end{corl}}
\def\ba{\begin{algo}}
\def\ea{\end{algo}}
\def\br{\begin{Remark}}
\def\er{\end{Remark}}
\def\bs{\begin{assum}}
\def\es{\end{assum}}
\def\bt{\begin{theorem}}
\def\et{\end{theorem}\vskip 3pt}
\def\bl{\begin{lemma}}
\def\el{\end{lemma}}
\def\ep{\end{proposition}}
\def\bp{\begin{proposition}}
\def\qed{\hfill{$\Box$}\vskip 5pt}
\def\be{\begin{example}}
\def\ee{\end{example}}
\def\bd{\begin{definition}}
\def\ed{\end{definition}}
\def\bc{\begin{corollary}}
\def\ec{\end{corollary}}
\def\proof{\noindent\it Proof. \hspace{1mm}\rm}
\input epsf
\begin{document}
\title{\bf Centrosymmetric, Skew Centrosymmetric and Centrosymmetric Cauchy Tensors}
\author{Haibin Chen\thanks{Department of Applied Mathematics, The Hong Kong Polytechnic University, Hung Hom,
Kowloon, Hong Kong. Email: chenhaibin508@163.com.}, \quad
Zhen Chen\thanks{School of Mathematics and Computer Science, Guizhou Normal
University, Guiyang, 550001, P. R. of China. Email: zchen@gznu.edu.cn. This author is supported by National Natural Science
Foundation of China ( Grant No.11201092, 11261012). His work was partially done when he was visiting The Hong Kong Polytechnic University}, \quad
Liqun Qi
\thanks{Department of Applied Mathematics, The Hong Kong Polytechnic University, Hung Hom,
Kowloon, Hong Kong. Email: maqilq@polyu.edu.hk. This author's work
was supported by the Hong Kong Research Grant Council (Grant No.
PolyU 502510, 502111, 501212 and 501913).} }
\maketitle
\begin{abstract} Recently, Zhao and Yang introduced centrosymmetric
tensors. In this paper, we further introduce skew centrosymmetric
tensors and centrosymmetric Cauchy tensors, and discuss properties
of these three classes of structured tensors.   Some sufficient and
necessary conditions for a tensor to be centrosymmetric or
skew centrosymmetric are given. We show that, a general
tensor can always be expressed as the sum of a centrosymmetric
tensor and a skew centrosymmetric tensor.   Some sufficient and
necessary conditions for a Cauchy tensor to be centrosymmetric
or skew centrosymmetric are also given. Spectral
properties on H-eigenvalues and H-eigenvectors of centrosymmetric, skew centrosymmetric and centrosymmetric Cauchy
tensors are discussed.   Some further questions on these tensors are
raised.


\end{abstract}

\newpage
\section{Introduction}

Let $\mathbb{R}^n$ be the $n$ dimensional real Euclidean space.
Denote the set of all natural numbers by $N$. Suppose $m$ and $n$
are two positive natural numbers and denote $[n]=\{1,2,\cdots,n\}$.

Centrosymmetric and skew centrosymmetric matrices play an
important role in information theory, linear system theory and
numerical analysis \cite{Andrew73, Anderson73,
Cant76,Demm97,Weaver85}. Discussion on various properties of such
matrices can be traced back to Muir \cite{Muir83}. Motivated by
these notions, Zhao and Yang introduced centrosymmetric tensors and
discussed properties of spectral radii of nonnegative
centrosymmetric tensors \cite{Zhao14}.

We now define centrosymmetric tensors and skew centrosymmetric
tensors.   The definition of centrosymmetric tensors here is the
same as Definition 2.1 of \cite{Zhao14}.

\bd\label{def11} Suppose an order $m$ dimension $n$ real tensor
$\mathcal{A}=(a_{i_1i_2\cdots i_m})$ satisfies
$$a_{i_1i_2\cdots i_m}=a_{n-i_1+1 n-i_2+1 \cdots n-i_m+1},~ i_j\in [n],~j\in [m].$$
Then $\mathcal{A}$ is called a centrosymmetric tensor. $\mathcal{A}$
is called a skew centrosymmetric tensor if it satisfies
$$a_{i_1i_2\cdots i_m}=-a_{n-i_1+1 n-i_2+1 \cdots n-i_m+1},~ i_j\in [n],~j\in [m].$$
\ed

By Definition \ref{def11}, a centrosymmetric tensor is symmetric
about its center. When dimension $n$ is odd, the centrosymmetric
tensor has the central entry $a_{ii\cdots i}$, where
$i=\frac{n+1}{2}$. When $n$ is even, there is no central entry. For
cases $m=2$, $n=2$ and $n=3$ respectively, we have
$$A=\left(
\begin{array}{cc}
a & b\\
b & a\\
\end{array}
\right),~~
B=\left(
\begin{array}{ccc}
a & b & c\\
d & e & d\\
c & b & a
\end{array}
\right).
$$
As we look at centrosymmetric tensors, we will find that they have
many interesting properties, comparable in some ways with symmetric
tensors. In general, a centrosymmetric tensor is not a symmetric
tensor. We have to point out that, in this paper, we always consider
order $m$ dimension $n$ centrosymmetric and skew centrosymmetric
tensors defined in the real field $\mathbb{R}$.

Apparently, centrosymmetric and skew centrosymmetric tensors are
structured tensors.   Recently, many interesting and impressed
properties of structured tensors have been discovered, and a lot of
research papers on structured tensors appeared
\cite{chen14,Chen14,Ding13,Ding14,HH14,Qi14,qi14,QXX,
Song14,song14,YY14,Zhang12}.    These include M tensors, circulant
tensors, completely positive tensors, Hankel tensors, Hilbert
tensors, P tensors, B tensors and Cauchy tensors. These papers not
only established results on spectral properties, positive
semi-definiteness and definiteness of structured tensors, but also
gave some important applications of structured tensors in stochastic
process and data fitting \cite{Chen14,Ding14}.

Our paper is organized as follows. In the next section, definitions
of tensor products, H-eigenvalues and
H-eigenvectors are given. In Section 3, basic properties of
centrosymmetric and skew centrosymmetric tensors are presented.
Firstly, we prove that the product of two arbitrary centrosymmetric
tensors is also centrosymmetric. Secondly, several
sufficient and necessary conditions for a tensor to be
centrosymmetric or skew centrosymmetric are given.
They are natural extensions of the matrix case. Furthermore, we
show that any general tensor can be denoted as the sum of a
centrosymmetric tensor and a skew centrosymmetric tensor. Some
properties on left inverses and right inverses of centrosymmetric and
skew centrosymmetric tensors are also presented in that section.
Properties on H-eigenvectors of centrosymmetric and skew
centrosymmetric tensors are discussed in Section 4. We prove that
some real lower dimensional tensors always have symmetric
H-eigenvectors or skew H-eigenvectors. It is proven that all
H-eigenvectors of a centrosymmetric tensor are still H-eigenvectors
of the tensor which is resulted from reversing the orders of the
entries. For a skew centrosymmetric tensor, all nonzero
H-eigenvalues must exist as pairs, which means that the reversed
value of a nonzero H-eigenvalue remains as an H-eigenvalue of that
tensor. In Section 5, the notion of centrosymmetric Cauchy tensor is
introduced.    It is proved that a Cauchy tensor is centrosymmetric
if and only if its generating vector is symmetric. We prove that
there is no odd dimension skew centrosymmetric Cauchy tensors.
Furthermore, when a centrosymmetric Cauchy tensor is of even order,
then its H-eigenvectors corresponding to any nonzero H-eigenvalues
are symmetric vectors. For a centrosymmetric Cauchy tensor of odd
order, the absolute vectors of H-eigenvectors corresponding to any
nonzero H-eigenvalues are symmetric. We conclude this paper with
some final remarks in Section 6.

By the end of the introduction, we add some comments on the notation
that will be used in the sequel. Let $\mathbb{C}^n$ be the $n$
dimensional complex space and let $\mathbb{C}$ be the complex field.
Vectors are denoted by italic lowercase letters i.e. $x,~ y,\cdots$,
and tensors are written as calligraphic capitals such as
$\mathcal{A}, \mathcal{T}, \cdots.$ Suppose $e\in \mathbb{R}^n$ be
all one vectors. Let $\mathcal{I}=(\delta_{i_1i_2\cdots i_m})$
denote the real identity tensor. If the symbol $|\cdot|$ is used on
a vector $x=(x_1,x_2,\cdots,x_n)$, then we get another vector
$|x|=(|x_1|,|x_2|,\cdots,|x_n|)$.

\setcounter{equation}{0}
\section{Preliminaries}
In this section, we present some basic definitions that will be used
in the sequel, such as tensor product, H-eigenvalue and
H-eigenvector.

\bd\label{def21}$^{\cite{Bu14}}$ Let $\mathcal{A}\in \mathbb{C}^{n_1\times n_2 \times \cdots \times n_2}$ and
$\mathcal{B}\in \mathbb{C}^{n_2\times n_3 \times \cdots \times n_{k+1}}$ be order $m\geq2$ and $k\geq1$ tensors,
respectively. The product $\mathcal{A}\mathcal{B}$ is the following tensor $\mathcal{C}$ of order $(m-1)(k-1)+1$
with entries:
$$c_{i\alpha_1\alpha_2\cdots\alpha_{m-1}}=\sum_{i_2,\cdots,i_m\in [n_2]}a_{ii_2\cdots i_m}b_{i_2\alpha_1}\cdots b_{i_m\alpha_{m-1}},$$
where $i\in [n_1], \alpha_1,\alpha_2,\cdots,\alpha_{m-1} \in [n_3]\times\cdots\times[n_{k+1}]$.
\ed

In this paper, we mainly study the case when $n_1=n_2=\cdots=n_{k+1}=n$. The product $\mathcal{A}\mathcal{B}$ was defined in \cite{Shao13,Shao132} when $n_1=n_2=\cdots=n_{k+1}=n$.

The definition of eigenvalue-eigenvector pairs of real symmetric
tensors comes from \cite{Qi05}.  Here we allow the tensors to be not
symmetric.

\bd\label{def22} Let $\mathbb{C}$ be the complex field. A pair
$(\lambda, x)\in \mathbb{C}\times \mathbb{C}^n\setminus \{0\}$ is
called an eigenvalue-eigenvector pair of a real tensor $\mathcal{T}$
with order $m$ dimension $n$, if they satisfy
\begin{equation}\label{e21}
\mathcal{T}x^{m-1}=\lambda x^{[m-1]},
\end{equation}
where $\mathcal{T}x^{m-1}=\left(\sum_{i_2,\cdots,i_m=1}^n t_{ii_2\cdots i_m}x_{i_2}\cdots x_{i_m} \right)_{1\leq i\leq n}$ and $x^{[m-1]}=(x_i^{m-1})_{1\leq i\leq n}$ are dimension $n$ vectors.
\ed
In Definition \ref{def22}, if $\lambda\in \mathbb{R}$ and the corresponding eigenvector $x\in \mathbb{R}^n$, then $\lambda,~x$ are called H-eigenvalue and
H-eigenvector respectively.

\setcounter{equation}{0}
\section{Basic Properties of Centrosymmetric and Skew Centrosymmetric Tensors}

In this section, we first give some results about products of
centrosymmetric tensors and skew centrosymmetric tensors. Then, some
sufficient and necessary conditions for a tensor to be a
centrosymmetric tensor or a skew centrosymmetric tensor are
presented, which are natural extensions of the matrix case. Finally,
we present properties of left inverses and right inverses of
centrosymmetric and skew centrosymmetric tensors.

\bl\label{lema31} Assume $B$ is an $n\times n$ square
centrosymmetric matrix and $\mathcal{A}$ is an order $m$ dimension
$n$ centrosymmetric tensor. Then $B\mathcal{A}$ is an order $m$
dimension $n$ centrosymmetric tensor. \el \proof By Definition
\ref{def21}, we have
$$(B\mathcal{A})_{i_1i_2\cdots i_m}=\sum_{j\in [n]}b_{i_1j}a_{ji_2\cdots i_m}.$$
For any $i_1,i_2,\cdots,i_m\in [n]$, since $B$ and $\mathcal{A}$ are centrosymmetric, so
$$\begin{array}{rl} (B\mathcal{A})_{i_1i_2\cdots i_m}=&\sum_{j\in [n]}b_{i_1j}a_{ji_2\cdots i_m} \\
=&\sum_{j\in [n]}b_{n-i_1+1 n-j+1}a_{n-j+1 n-i_2+1 \cdots n-i_m+1}\\
=&\sum_{l\in [n]}b_{n-i_1+1 l}a_{l n-i_2+1 \cdots n-i_m+1}\\
=&(B\mathcal{A})_{n-i_1+1 n-i_2+1 \cdots n-i_m+1}.
\end{array}$$
Combining this with Definition \ref{def11}, we know that $B\mathcal{A}$ is a centrosymmetric tensor.
\qed

\bl\label{lema32} Suppose $\mathcal{A}$ and $B$ are defined as in
Lemma \ref{lema31}. Then $\mathcal{A}B$ is a centrosymmetric tensor.
\el \proof By Definition \ref{def21} and the fact that $\mathcal{A}$
and $B$ are centrosymmetric, we have
$$\begin{array}{rl} (\mathcal{A}B)_{i_1i_2\cdots i_m}=&\sum_{j_2,j_3,\cdots, j_m\in [n]} a_{i_1j_2\cdots j_m}b_{j_2i_2}\cdots b_{j_mi_m}\\
=&\sum_{j_2,j_3,\cdots, j_m\in [n]}a_{n-i_1+1 n-j_2+1 \cdots n-j_m+1}b_{n-j_2+1 n-i_2+1}\cdots b_{n-j_m+1 n-i_m+1}\\
=&\sum_{l_2,l_3,\cdots, l_m\in [n]}a_{n-i_1+1 l_2 \cdots l_m}b_{l_2 n-i_2+1}\cdots b_{l_m n-i_m+1}\\
=&(\mathcal{A}B)_{n-i_1+1 n-i_2+1 \cdots n-i_m+1},
\end{array}
$$
for any $i_1,i_2,\cdots,i_m\in [n]$. Thus $\mathcal{A}B$ is a
centrosymmetric tensor. \qed \bt\label{thm31} Let  $\mathcal{A}$ be
order $m$ dimension $n$ tensor and $\mathcal{B}$ be order $k$
dimension $n$ tensor. Assume $\mathcal{A}$ and $\mathcal{B}$ are
centrosymmetric tensors. Then the production $\mathcal{A}
\mathcal{B}$ is an order $(m-1)(k-1)+1$ dimension $n$
centrosymmetric tensor. \et \proof By Definition \ref{def21}, for
any $i_1\in [n]$, $\alpha_j=\alpha^j_1 \alpha^j_2 \cdots
\alpha^j_{k-1} \in [n]^{k-1},~j\in [m-1]$, we have
$$\begin{array}{rl} (\mathcal{A}\mathcal{B})_{i_1\alpha_1 \cdots \alpha_{m-1}}=&\sum_{j_2,j_3,\cdots, j_m\in [n]} a_{i_1j_2\cdots j_m}b_{j_2\alpha_1}\cdots b_{j_m\alpha_{m-1}}\\
=&\sum_{j_2,j_3,\cdots, j_m\in [n]}a_{n-i_1+1 n-j_2+1 \cdots n-j_m+1}b_{n-j_2+1 n-\alpha_1+1}\cdots b_{n-j_m+1 n-\alpha_{m-1}+1}\\
=&\sum_{l_2,l_3,\cdots, l_m\in [n]}a_{n-i_1+1 l_2 \cdots l_m}b_{l_2 n-\alpha_1+1}\cdots b_{l_m n-\alpha_{m-1}+1}\\
=&(\mathcal{A}B)_{n-i_1+1 n-\alpha_1+1 \cdots n-\alpha_{m-1}+1},
\end{array}
$$
where $n-\alpha_j+1$ means $n-\alpha^j_t+1$ for every index $\alpha^j_t$ in $\alpha_j$,~$t\in [k-1]$. Here, the second equality follows that
$\mathcal{A}$ and $\mathcal{B}$ are centrosymmetric tensors. Obviously $\mathcal{A}\mathcal{B}$ are centrosymmetric tensors.
\qed

From the proof process of Theorem \ref{thm31}, we have the following
corollaries and we omit the proofs for simplicity.

\bc\label{corol31} Suppose tensor $\mathcal{A}$, $\mathcal{B}$ are defined as in Theorem \ref{thm31}. Then the following statements holds:

(i) if $\mathcal{A}$ is skew centrosymmetric and $\mathcal{B}$ is centrosymmetric, then $\mathcal{A}\mathcal{B}$ is skew centrosymmetric.

(ii) if $\mathcal{A}$ is centrosymmetric and $\mathcal{B}$ is skew centrosymmetric, then $\mathcal{A}\mathcal{B}$ is centrosymmetric when $m$ is odd; $\mathcal{A}\mathcal{B}$ is skew centrosymmetric when $m$ is even.

(iii) if $\mathcal{A}$ and $\mathcal{B}$ are both skew ceontrosymmetric, then $\mathcal{A}\mathcal{B}$ is centrosymmetric when $m$ is even; $\mathcal{A}\mathcal{B}$ is skew centrosymmetric when $m$ is odd.
\ec

\bc\label{corol32} For any finite dimension $n$ tensors
$\mathcal{A}_1,~\mathcal{A}_2,\cdots,~\mathcal{A}_s$, if they are
all centrosymmetric tensors, then the product
$\mathcal{A}_1\mathcal{A}_2\cdots\mathcal{A}_s$ is also a
centrosymmetric tensor. \ec

Let $r_i,~i\in [n]$ denote the sum of some elements in $\mathcal{A}$
such that
$$r_i=\sum_{i_2,i_3,\cdots,i_m \in [n]}a_{ii_2\cdots i_m},~i\in [n].$$
By the definition of centrosymmetric tensors and skew
centrosymmetric tensors, we have the following conclusions.

\bt\label{thm32} Let $\mathcal{A}=(a_{i_1i_2\cdots i_m})$ be an
order $m$ dimension $n$ tensor. If $A$ is centrosymmetric, then
$r_i=r_{n-i+1}$; if $A$ is skew centrosymmetric, then
$r_i=-r_{n-i+1}$.\et

\bc\label{corl33} Assume $\mathcal{A}$ is defined as in Theorem
\ref{thm32}. Suppose $\mathcal{A}$ is skew centrosymmetric and $n$
is an odd number.  Then there are at least one zero element in
$\mathcal{A}$ and at least one $i\in [n]$ satisfying $r_i=0$. \ec
\proof From Definition \ref{def11} and the fact that $n$ is odd, let
$i=\frac{n+1}{2}$, then we have
$$a_{ii\cdots i}=-a_{ii\cdots i},~r_i=-r_i,$$
which implies that
$$a_{ii\cdots i}=0,~r_i=0$$
and the desired results hold.
\qed

We now give some sufficient and necessary conditions for a tensor to
be centrosymmetric or skew centrosymmetric. Let $J$ be the $n\times
n$ real matrix with elements satisfying $J_{ij}=\delta_{i
n-j+1},~1\leq i, j\leq n,$ where $\delta_{i n-j+1}$ denotes the
Kronecker delta
$$
J=\left(
\begin{array}{ccccc}
0 & 0 & \cdots & 0 & 1\\
0 & 0 & \cdots & 1 & 0\\
  &   & \cdots &   &  \\
0 & 1 & \cdots & 0 & 0\\
1 & 0 & \cdots & 0 & 0\\
\end{array}
\right).
$$

\bt\label{thm33} Let $\mathcal{A}$ be an order $m$ dimension $n$ tensor. Then $\mathcal{A}$ is centrosymmetric
if and only if $J\mathcal{A}J=\mathcal{A}$; $\mathcal{A}$ is skew centrosymmetric
if and only if $J\mathcal{A}J=-\mathcal{A}$.
\et

\proof For any $i_j\in [n],~j\in [m]$, by Definition \ref{def21}, we have
$$
(J\mathcal{A}J)_{i_1 i_2 \cdots i_m}=\sum_{j_1,j_2,\cdots,j_m \in [n]}J_{i_1j_1}a_{j_1 j_2 \cdots j_m}J_{j_2i_2}\cdots J_{j_mi_m}
$$
By Definition \ref{def11} and the definition of matrix $J$, one has
$$
\begin{array}{rl}
(J\mathcal{A}J)_{i_1 i_2 \cdots i_m}=&\sum_{j_1,j_2,\cdots,j_m \in [n]}J_{i_1j_1}a_{j_1 j_2 \cdots j_m}J_{j_2i_2}\cdots J_{j_mi_m}\\
=&a_{n-i_1+1 n-i_2+1 \cdots n-i_m+1},
\end{array}
$$
which implies that the sufficient and necessary condition holds. Moreover, the second conclusion can be proven similarly.
\qed

Since $JJ=I$, where $I$ is the $n\times n$ identity matrix,
according to Proposition 1.1 of \cite{Shao13} and Theorem 1.1 of
\cite{Shao13}, we have the following conclusion.

\bt\label{thm34} Let $\mathcal{A}$ be an order $m$ dimension $n$ tensor. Then $\mathcal{A}$ is centrosymmetric
if and only if $\mathcal{A}J=J\mathcal{A}$; $\mathcal{A}$ is skew centrosymmetric
if and only if $\mathcal{A}J=-J\mathcal{A}$.
\et

Let $x=(x_1,x_2,\cdots,x_n)\in \mathbb{R}^n$. Then $Jx$ is a vector that can be gotten by reversing orders of elements of $x$.  If $Jx=x$, we call $x$ is a symmetric vector and it is called skew symmetric if
$Jx=-x$. For any given tensor $\mathcal{A}=(a_{i_1i_2\cdots i_m})$ with order $m$ dimension $n$, the corresponding homogeneous polynomial is denoted by
$$f(x)=\mathcal{A}x^m=\sum_{i_1,i_2,\cdots, i_m\in [n]}a_{i_1i_2\cdots i_m}x_{i_1}\cdots x_{i_m}.$$

\bt\label{thm35} Suppose order $m$ dimension $n$ tensor
$\mathcal{A}$ is centrosymmetric. Then $f(Jx)=f(x)$ for any $x\in
\mathbb{R}^n$; If $\mathcal{A}$ is skew centrosymmetric, then
$f(Jx)=-f(x)$. \et \proof Let $y=Jx=(x_n,x_{n-1},\cdots,x_2,x_1)$,
which means $y_i=x_{n-i+1},~i\in [n]$. If $\mathcal{A}$ is
centrosymmetric, then we have
\begin{equation}\label{e31}
\begin{array}{rl} f(Jx)=&f(y)\\
=&\sum_{i_1,i_2,\cdots, i_m\in [n]}a_{i_1i_2\cdots i_m}y_{i_1}\cdots y_{i_m}\\
=&\sum_{i_1,i_2,\cdots, i_m\in [n]}a_{i_1i_2\cdots i_m}x_{n-i_1+1}x_{n-i_2+1}\cdots x_{n-i_m+1}\\
=&\sum_{i_1,i_2,\cdots, i_m\in [n]}a_{n-i_1+1 n-i_2+1 \cdots n-i_m+1}x_{n-i_1+1}x_{n-i_2+1}\cdots x_{n-i_m+1}\\
=&\sum_{j_1,j_2,\cdots, j_m\in [n]}a_{j_1 j_2 \cdots j_m}x_{j_1}x_{j_2}\cdots x_{j_m}\\
=&f(x).
\end{array}
\end{equation}
When $\mathcal{A}$ is skew centrosymmetric, one has
\begin{equation}\label{e32}
\begin{array}{rl} f(Jx)=&f(y)\\
=&\sum_{i_1,i_2,\cdots, i_m\in [n]}a_{i_1i_2\cdots i_m}y_{i_1}\cdots y_{i_m}\\
=&\sum_{i_1,i_2,\cdots, i_m\in [n]}a_{i_1i_2\cdots i_m}x_{n-i_1+1}x_{n-i_2+1}\cdots x_{n-i_m+1}\\
=&-\sum_{i_1,i_2,\cdots, i_m\in [n]}a_{n-i_1+1 n-i_2+1 \cdots n-i_m+1}x_{n-i_1+1}x_{n-i_2+1}\cdots x_{n-i_m+1}\\
=&-\sum_{j_1,j_2,\cdots, j_m\in [n]}a_{j_1 j_2 \cdots j_m}x_{j_1}x_{j_2}\cdots x_{j_m}\\
=&-f(x).
\end{array}
\end{equation}
By (\ref{e31}) and (\ref{e32}), we know that the desired results hold. \qed

Suppose $\mathcal{A}=(a_{i_1i_2\cdots i_m})$ and $\mathcal{B}=(b_{i_1i_2\cdots i_m})$ are two order $m$ dimension $n$ tensors, the Hadamard product
of $\mathcal{A}$ and $\mathcal{B}$ is defined as
\begin{equation}\label{e33}
\mathcal{A}\circ \mathcal{B}=(a_{i_1i_2\cdots i_m}b_{i_1i_2\cdots i_m}),
\end{equation}
which is still an order $m$ dimension $n$ tensor. Now, we present
several conclusions about the Hadamard product of centrosymmetric
tensors and skew centrosymmetric tensors.

\bt\label{thm36} For two order $m$ dimension $n$ tensors
$\mathcal{A}$ and $\mathcal{B}$, we have the following statements:

(i) if $\mathcal{A}$ and $\mathcal{B}$ are centrosymmetric, then
$\mathcal{A}\circ \mathcal{B}$ is centrosymmetric;

(ii) if $\mathcal{A}$ and $\mathcal{B}$ are skew centrosymmetric
tensors, then $\mathcal{A}\circ \mathcal{B}$ is centrosymmetric;

(iii) if $\mathcal{A}$ is centrosymmetric  and $\mathcal{B}$ is skew
centrosymmetric, then $\mathcal{A}\circ \mathcal{B}$ is skew
centrosymmetric. \et

\proof By Definition \ref{def11} and (\ref{e33}), it is easy to check the authenticity of the results. Thus, we omit the proof.\qed

As we all know that any matrix can be decomposed to the sum of a symmetric matrix and a skew symmetric matrix. Similarly, we have the following result.

\bt\label{thm37} Any order $m$ dimension $n$ tensor $\mathcal{A}$
can be expressed as the sum of a centrosymmetric tensor and a skew
centrosymmetric tensor. \et

\proof Without loss of generality, let $\mathcal{A}=(a_{i_1 i_2 \cdots i_m})$, $i_j\in [n],~j\in [m]$. Set a new tensor $\mathcal{A}^c=(a^c_{i_1 i_2 \cdots i_m})$ such that
$$a^c_{i_1 i_2 \cdots i_m}=a_{n-i_1+1 n-i_2+1 \cdots n-i_m+1},~i_j\in [n],~j\in [m].$$
From a direct computation, we have
$$\mathcal{A}=\frac{\mathcal{A}+\mathcal{A}^c}{2}+\frac{\mathcal{A}-\mathcal{A}^c}{2},$$
where $\frac{\mathcal{A}+\mathcal{A}^c}{2}$ is centrosymmetric and $\frac{\mathcal{A}-\mathcal{A}^c}{2}$ is skew centrosymmetric. Thus, the desired result follows.\qed

Another important property of centrosymmetric matrices is that the
inverse matrix of a centrosymmetric matrix is also centrosymetric
\cite{Weaver85}. So, we want to know whether the inverse of a
centrosymmetric tensor is centrosymmetric or not. Unfortunately,
there is no definition of the inverse of a tensor. But, definitions
of left inverse and right inverse of tensors are given in
\cite{Bu14}. In the following, we will study the centrosymmetric
property of left inverse tensors and right inverse tensors under the
assumption that a centrosymmetric tensor has left inverse and right
inverse.

In \cite{Bu14}, Bu C. et al. presented the definition of left inverse and right inverse of tensors as below.

\bd\label{def31}$^{\cite{Bu14}}$ Let $\mathcal{A}$ be a tensor of order $m$ and dimension $n$ and let $\mathcal{B}$
be a tensor of order $k$ and dimension $n$. If $\mathcal{A}\mathcal{B}=\mathcal{I}$, then $\mathcal{A}$ is called an order $m$ left inverse of
$\mathcal{B}$, and $\mathcal{B}$ is called an order $k$ right inverse of $\mathcal{A}$.
\ed

\bt\label{thm38} Assume $\mathcal{A}=(a_{i_1i_2\cdots i_m})$ be a diagonal centrosymmetric tensor of order $m$ dimension $n$.
Then,

(i) $\mathcal{A}$ has real centrosymmetric left inverse if and only if $\mathcal{A}$ has nonzero diagonal entries;

(ii) when $m$ is even, $\mathcal{A}$ has real centrosymmetric right inverse if and only if $\mathcal{A}$  has nonzero diagonal entries;

(iii) when $m$ is odd, $\mathcal{A}$ has real centrosymmetric right inverse if all diagonal entries of $\mathcal{A}$ are positive.
\et
\proof (i) By Definition \ref{def31}, $\mathcal{A}$ has real centrosymmetric left inverse if and only if there exists a real centrosymmetric
tensor $\mathcal{B}=(b_{i_1i_2\cdots i_k})$ with order $k$ and dimension $n$ such that
\begin{equation}\label{e34}
\mathcal{B}\mathcal{A}=\mathcal{I},
\end{equation}
and
$$b_{i_1i_2\cdots i_k}=b_{n-i_1+1 n-i_2+1 \cdots n-i_m+1}.$$
For any $i\in [n]$, $\alpha_j \in [n]^{m-1},~j\in [k-1]$we have
$$
\begin{array}{rl} (\mathcal{B}\mathcal{A})_{i\alpha_1 \cdots \alpha_{k-1}}=&\sum_{j_2,j_3,\cdots,j_k\in [n]}b_{ij_2\cdots j_k}a_{j_2\alpha_1}\cdots a_{j_k\alpha_{k-1}}\\
=&\delta_{i\alpha_1 \cdots \alpha_{k-1}}.
\end{array}$$
When $\alpha_j=ii\cdots i$ for all $j\in [k-1]$, one has
$$b_{ii\cdots i}a^{k-1}_{ii\cdots i}=1.$$
Thus, the existence of left inverse of $\mathcal{A}$ implies that all diagonal elements of $\mathcal{A}$ must be nonzero and the only if part holds.
For sufficient condition, if
$$a_{ii\cdots i}\neq0,~i\in [n],$$
let
$$b_{ii\cdots i}=\frac{1}{a^{k-1}_{ii\cdots i}},~i\in [n]$$
and $b_{i_1i_2\cdots i_k}=0$ for the others. Then, $\mathcal{B}$ is centrosymmetric since $\mathcal{A}$ is centrosymmetric and it is easy to check equation
(\ref{e34}) holds. Thus tensor $\mathcal{B}$ is an order $k$ real left inverse of $\mathcal{A}$.

(ii) For only if part, there is an order $k$ dimension $n$ real centrosymmetric tensor $\mathcal{B}=(b_{i_1i_2\cdots i_k})$ such that
$$
\begin{array}{rl}(\mathcal{A}\mathcal{B})_{i\alpha_1 \cdots \alpha_{m-1}}=&\sum_{j_2,j_3,\cdots, j_m\in [n]} a_{ij_2\cdots j_m}b_{j_2\alpha_1}\cdots b_{j_m\alpha_{m-1}}\\
=&\delta_{i\alpha_1 \cdots \alpha_{m-1}}.
\end{array}$$
for $i\in [n]$, $\alpha_j \in [n]^{k-1},~j\in [m-1]$. For diagonal entries of $\mathcal{A}\mathcal{B}$, we have
\begin{equation}\label{e35}
a_{ii\cdots i}b^{m-1}_{ii\cdots i}=\delta_{ii\cdots i}=1,~i\in [n],
\end{equation}
which implies that tensor $\mathcal{A}$ has nonzero diagonal entries.

For sufficient conditions, let the entries of tensor $\mathcal{B}$ be that
$$b_{ii\cdots i}=(\frac{1}{a_{ii\cdots i}})^{\frac{1}{m-1}}, ~i\in [n]$$
and $b_{i_1i_2\cdots i_k}=0$ for the others. Then, by a direct computation, we know that $\mathcal{B}$ is a real centrosymmetric right inverse of $\mathcal{A}$.

(iii) When $m$ is odd, by (\ref{e35})
$$a_{ii\cdots i}b^{m-1}_{ii\cdots i}=\delta_{ii\cdots i}=1,~i\in [n],$$
we have that all diagonal elements of $\mathcal{A}$ are positive. The others are similar to the proof of (ii).
\qed

\bt\label{thm39} Suppose $\mathcal{A}$ is a centrosymmetric tensor of order $m$ and dimension $n$. If $\mathcal{A}$ has an order 2 dimension $n$
real left inverse, then it must be unique and centrosymmetric.
\et

\proof Suppose matrix $B$ is an order 2 real left inverse of tensor $\mathcal{A}$. By Definition \ref{def31}, we have
$$B\mathcal{A}=\mathcal{I}.$$
From Proposition 2.1 of \cite{Shao13} and Problem 1 of \cite{Shao132}, we obtain
$$det(B)\neq0,$$
which means that $B$ is a nonsingular matrix. Let $B^{-1}=(b^{-1}_{ij})$ denote the inverse of $B$. From Theorem 1.1 of \cite{Shao13}, one has
$$\mathcal{A}=B^{-1}\mathcal{I}.$$
Thus, for any $i,j\in [n]$, it holds that
$$
a_{ijj\cdots j}=\sum_{t\in [n]}b^{-1}_{it}\delta_{tjj\cdots j}=b^{-1}_{ij}
$$
and
$$
a_{n-i+1 n-j+1\cdots n-j+1}=\sum_{t\in [n]}b^{-1}_{n-i+1 t}\delta_{t n-j+1 n-j+1 \cdots n-j+1}=b^{-1}_{n-i+1 n-j+1}.
$$
Since tensor $\mathcal{A}$ is centrosymmetric, so
$$b^{-1}_{n-i+1 n-j+1}=a_{n-i+1 n-j+1\cdots n-j+1}=a_{ijj\cdots j}=b^{-1}_{ij},$$
which implies that $B^{-1}$ is a centrosymmetric matrix. By Proposition 6 of \cite{Weaver85}, we know that $B$ is centrosymmetric.

Assume $\mathcal{A}$ has another order 2 real left inverse $C$. Then,
$$\mathcal{A}=B^{-1}\mathcal{I}=C^{-1}\mathcal{I},$$
where $C^{-1}$ is the inverse of $C$. Then,
$$(B^{-1}-C^{-1})\mathcal{I}=0.$$
Combining this with Lemma 2.1 of \cite{Bu14}, we have
$$B^{-1}=C^{-1}.$$
By the fact that a nonsingular matrix has a unique inverse matrix, it follows that $B=C$ and the desired results hold.
\qed

\bt\label{thm310}  Suppose $\mathcal{A}$ is a centrosymmetric tensor of order $m$ and dimension $n$. Let $m$ be even. If $\mathcal{A}$ has an order 2 dimension $n$
real right inverse, then it must be unique and centrosymmetric.
\et
\proof Let $B=(b_{ij})$ be any order 2 real right inverse of $\mathcal{A}$. By Proposition 2.1 of \cite{Shao13} and Problem 1 of \cite{Shao132}, we
know that $B$ is nonsingular. So, from Theorem 1.1 of \cite{Shao13}, we obtain
$$\mathcal{A}B=\mathcal{I},$$
which can be written
$$\mathcal{A}=\mathcal{I}B^{-1},$$
where $B^{-1}=(b^{-1}_{ij})$ is the inverse of matrix $B$.
For any $i,j\in [n]$, one has
$$a_{ijj\cdots j}=\sum_{i_2,i_3,\cdots,i_m \in [n]}\delta_{ii_2\cdots i_m}b^{-1}_{i_2j}b^{-1}_{i_3j}\cdots b^{-1}_{i_mj}.$$
Thus, we obtain
$$a_{ijj\cdots j}=( b^{-1}_{ij})^{m-1},~a_{n-i+1 n-j+1 n-j+1\cdots n-j+1}=( b^{-1}_{n-i+1 n-j+1})^{m-1}.$$
By the fact that tensor $\mathcal{A}$ is centrosymmetric, it follows that
$$( b^{-1}_{ij})^{m-1}=( b^{-1}_{n-i+1 n-j+1})^{m-1},~i,j\in [n],$$
and
$$ b^{-1}_{ij}= b^{-1}_{n-i+1 n-j+1},~i,j\in [n],$$
since $m$ is even. So matrix $B^{-1}$ is centrosymmetric and $B$ is centrosymmetric from Proposition 6 of \cite{Weaver85}.

Assume $\mathcal{A}$ has another order 2 real right inverse $C$. Then,
$$\mathcal{A}=\mathcal{I}B^{-1}=\mathcal{I}C^{-1},$$
where $C^{-1}$ is the inverse of $C$. Then,
$$\mathcal{I}(B^{-1}-C^{-1})=0.$$
By Lemma 2.2 of \cite{Bu14}, we know that $B^{-1}=C^{-1}$ and $B=C$.
\qed

\setcounter{equation}{0}
\section{Spectral Properties of Centrosymmetric and Skew Centrosymmetric Tensors}

In this section, we present several conclusions about H-eigenvalues and H-eigenvectors of real centrosymmetric and skew
centrosymmetric tensors.

In \cite{Weaver85}, it listed that all H-eigenvectors of real matrices with dimension $2\times2$ or dimension $3\times 3$ are either symmetric or
skew symmetric. But, these nice formulas cannot be extended to the $4\times 4$ case. Next, we will give some properties about H-eigenvectors
of dimension 2 and dimension 3 centrosymmetric tensors. The following two theorems show that order $m$ dimension 2 and order $m$ dimension 3 centrosymmetric
tensors always have symmetric H-eigenvectors or skew symmetric H-eigenvectors respectively.

\bt\label{thm41} Suppose $\mathcal{A}=(a_{i_1i_2\cdots i_m})$ is a centrosymmetric tensor of order m and dimension 2. Then, $\sum\limits_{i_2,\cdots,i_m=1}^{2} {a}_{1i_2 \cdots i_m}$ and $\sum\limits_{i_2,\cdots,i_m=1}^{2}{a}_{1i_2 \cdots i_m}(-1)^{i_2+\cdots+i_m+m-1}$ are H-eigenvalues of $\mathcal{A}$ with symmetric
H-eigenvector and skew symmetric H-eigenvector respectively.
\et
\proof Let $e=(1,1)^T$ and $u=(1,-1)^T$. From Definition \ref{def22} and the fact that $\mathcal{A}$ is centrosymmetric,
by a direct computation we have
$$\mathcal{A}e^{m-1}=(\sum\limits_{i_2,\cdots,i_m=1}^{2} {a}_{1i_2 \cdots i_m})e^{[m-1]}$$
and
$$\mathcal{A}u^{m-1}=(\sum\limits_{i_2,\cdots,i_m=1}^{2}{a}_{1i_2 \cdots i_m}(-1)^{i_2+\cdots+i_m+m-1})u^{[m-1]},$$
which imply that the desired results hold.
\qed

\bt\label{thm42} Assume $\mathcal{A}=(a_{i_1i_2\cdots i_m})$ is a centrosymmetric tensor of order m dimension 3. Suppose $m$ is even. Then, $\lambda=\sum\limits_{ i_2,\cdots,i_m\in \{1,3\}}a_{1i_2 \cdots i_m}(-1)^{r(i_2 \cdots i_m)}$ is an H-eigenvalue of $\mathcal{A}$ with skew symmetric H-eigenvector, where $r(i_2 \cdots i_m)$ denote the number of indices $i_2,\cdots,i_m$ equal 3.
\et
\proof Let $x=(1,0,-1)^T$. By Definition \ref{def22}, it is easy to check that
\begin{equation}\label{e41}
\begin{array}{rl} (\mathcal{A}x^{m-1})_1=&\sum\limits_{i_2,\cdots,i_m=1}^{3}a_{1i_2 \cdots i_m} x_{i_2} \cdots x_{i_m}\\
=&\sum\limits_{ i_2,\cdots,i_m\in \{1,3\}}a_{1i_2 \cdots i_m}x_{i_2} \cdots x_{i_m}\\
=&\sum\limits_{ i_2,\cdots,i_m\in \{1,3\}}a_{1i_2 \cdots i_m}(-1)^{r(i_2 \cdots i_m)}.
\end{array}
\end{equation}
Combining this with the fact $\mathcal{A}$ is centrosymmetric and $m$ is even, one has
\begin{equation}\label{e42}
(\mathcal{A}x^{m-1})_3=-(\mathcal{A}x^{m-1})_1.
\end{equation}
On the other hand,
$$
\begin{array}{rl} (\mathcal{A}x^{m-1})_2=&\sum\limits_{ i_2,\cdots,i_m\in \{1,3\}}a_{2i_2 \cdots i_m}x_{i_2} \cdots x_{i_m}\\
=&\sum\limits_{ i_2,\cdots,i_m\in \{1,3\}}a_{2i_2 \cdots i_m}(-1)^{r(i_2 \cdots i_m)}\\
=&\sum\limits_{ i_2,\cdots,i_m\in \{1,3\}}a_{2 4-i_2 \cdots 4-i_m}(-1)^{r(i_2 \cdots i_m)}\\
=&\sum\limits_{ j_2,\cdots,j_m\in \{1,3\}}a_{2 j_2 \cdots j_m}(-1)^{m-1-r(j_2 \cdots j_m)}\\
=&-\sum\limits_{ j_2,\cdots,j_m\in \{1,3\}}a_{2 j_2 \cdots j_m}(-1)^{r(j_2 \cdots j_m)}\\
=&-(\mathcal{A}x^{m-1})_2.
\end{array}
$$
Thus,
\begin{equation}\label{e43}
(\mathcal{A}x^{m-1})_2=0.
\end{equation}
By (\ref{e41})-(\ref{e43}), we know that $\lambda$ is an H-eigenvalue of $\mathcal{A}$ corresponding to the skew symmetric H-eigenvector $x$.
\qed

Now, we consider general order $m$ dimension $n$ centrosymmetric tensors and skew centrosymmetric tensors.
We will show that all H-eigenvectors of a centrosymmetric tensor are still H-eigenvectors
of the tensor which is resulted from reversing the orders of the
entries.   On the other side, for a skew centrosymmetric tensor, if it has a nonzero H-eigenvalue $\lambda$, then $-\lambda$ is still an H-eigenvalue of that
skew centrosymmetric tensor.

\bt\label{thm43} Let tensor $\mathcal{A}$ be a centrosymmetric tensor of order $m$ dimension $n$. If $\mathcal{A}$
has an H-eigenvalue $\lambda$ with an H-eigenvector $x$, then $Jx$ is also an H-eigenvector of $\mathcal{A}$ corresponding to $\lambda$.
\et
\proof By Definition \ref{def22}, we have
$$\mathcal{A}x^{m-1}=\lambda x^{[m-1]}.$$
Let $x=(x_1,x_2,\cdots,x_n)$, then $Jx=(x_n,x_{n-1},\cdots,x_1)$. For any $i\in [n]$, one has
\begin{equation}\label{e44}
\begin{array}{rl} (\mathcal{A}(Jx)^{m-1})_i=&\sum_{i_2,i_3,\cdots,i_m\in [n]}a_{ii_2i_3\cdots i_m}(Jx)_{i_2}\cdots (Jx)_{i_m}\\
=&\sum_{i_2,i_3,\cdots,i_m\in [n]}a_{ii_2i_3\cdots i_m}x_{n-i_2+1}x_{n-i_3+1}\cdots x_{n-i_m+1}\\
=&\sum_{i_2,i_3,\cdots,i_m\in [n]}a_{n-i+1 n-i_2+1\cdots n-i_m+1}x_{n-i_2+1}x_{n-i_3+1}\cdots x_{n-i_m+1}\\
=&\sum_{t_2,t_3,\cdots,t_m\in [n]}a_{n-i+1 t_2t_3\cdots t_m}x_{t_2}x_{t_3}\cdots x_{t_m}\\
=&\sum_{t_2,t_3,\cdots,t_m\in [n]}\lambda x_{n-i+1}^{m-1}\\
=&\sum_{t_2,t_3,\cdots,t_m\in [n]}\lambda (Jx)_i^{m-1}.
\end{array}
\end{equation}
Thus, $Jx$ is an H-eigenvector of $\mathcal{A}$ corresponding to the H-eigenvalue $\lambda$.
\qed

\bt\label{thm44} Let tensor $\mathcal{A}$ be a skew centrosymmetric tensor of order $m$ dimension $n$. If $\mathcal{A}$
has a nonzero H-eigenvalue $\lambda$ with an H-eigenvector $x$, then $Jx$ is also an H-eigenvector of $\mathcal{A}$ corresponding to the H-eigenvalue $-\lambda$.
\et
\proof By definition of H-eigenvalues and H-eigenvectors, we have
$$\mathcal{A}x^{m-1}=\lambda x^{[m-1]}.$$
Similarly, by (\ref{e44}), for any $i\in [n]$, one has
$$
\begin{array}{rl} (\mathcal{A}(Jx)^{m-1})_i=&\sum_{i_2,i_3,\cdots,i_m\in [n]}a_{ii_2i_3\cdots i_m}(Jx)_{i_2}\cdots (Jx)_{i_m}\\
=&\sum_{i_2,i_3,\cdots,i_m\in [n]}a_{ii_2i_3\cdots i_m}x_{n-i_2+1}x_{n-i_3+1}\cdots x_{n-i_m+1}\\
=&-\sum_{i_2,i_3,\cdots,i_m\in [n]}a_{n-i+1 n-i_2+1\cdots n-i_m+1}x_{n-i_2+1}x_{n-i_3+1}\cdots x_{n-i_m+1}\\
=&-\sum_{t_2,t_3,\cdots,t_m\in [n]}a_{n-i+1 t_2t_3\cdots t_m}x_{t_2}x_{t_3}\cdots x_{t_m}\\
=&-\sum_{t_2,t_3,\cdots,t_m\in [n]}\lambda x_{n-i+1}^{m-1}\\
=&\sum_{t_2,t_3,\cdots,t_m\in [n]}(-\lambda) (Jx)_i^{m-1}.
\end{array}
$$
Thus, $Jx$ is an H-eigenvector of $\mathcal{A}$ corresponding to the H-eigenvalue $-\lambda$.
\qed

\bt\label{thm45}
 Let $\mathcal{A}=(a_{i_1 i_2 \cdots i_m})$ be a centrosymmetric tensor of order $m$ dimension $n$. Then, all H-eigenvectors corresponding to the
H-eigenvalue $\lambda$, where $\dim \ker(\lambda \mathcal{I}-\mathcal{A})=1$, are either symmetric or skew-symmetric.
\et
\proof Suppose $x\in \mathbb{R}^n$ is a H-eigenvector of $\mathcal{A}$ corresponding to $\lambda$, where $\dim \ker(\lambda \mathcal{I}-\mathcal{A})=1$.
By Definition \ref{def22}, we have
\begin{equation}\label{e45}
\mathcal{A}x^{m-1}=\lambda x^{[m-1]}.
\end{equation}
By Theorem \ref{thm33} and (\ref{e45}), one has
$$
J\mathcal{A}Jx^{m-1}=\lambda x^{[m-1]},
$$
which implies that
\begin{equation}\label{e46}
\mathcal{A}Jx^{m-1}=\lambda Jx^{[m-1]}=\lambda (Jx)^{[m-1]}.
\end{equation}

On the other hand, by Definition \ref{def11}, we have
$$
\begin{array}{rl} (\mathcal{A}Jx^{m-1})_i=&\sum_{i_2,i_3,\cdots,i_m \in [n]}(\mathcal{A}J)_{ii_2i_3\cdots i_m}x_{i_2}x_{i_3}\cdots x_{i_m} \\
=&\sum_{i_2,i_3,\cdots,i_m \in [n]}(\sum_{j_2,j_3,\cdots,j_m \in [n]}a_{ij_2\cdots j_m}J_{j_2i_2}\cdots J_{j_mi_m})x_{i_2}x_{i_3}\cdots x_{i_m}\\
=&\sum_{i_2,i_3,\cdots,i_m \in [n]}a_{i n-i_2+1 \cdots n-i_m+1}x_{i_2}x_{i_3}\cdots x_{i_m}\\
=&\sum_{t_2,t_3,\cdots,t_m \in [n]}a_{i t_2t_3\cdots t_m}x_{n-t_2+1}x_{n-t_3+1}\cdots x_{n-t_m+1}\\
=&\sum_{t_2,t_3,\cdots,t_m \in [n]}a_{i t_2t_3\cdots t_m}(Jx)_{t_2}(Jx)_{t_3}\cdots (Jx)_{t_m}\\
=&(\mathcal{A}(Jx)^{m-1})_i,
\end{array}
$$
for any $i\in [n]$. So, it holds that
\begin{equation}\label{e47}
\mathcal{A}Jx^{m-1}=\mathcal{A}(Jx)^{m-1}.
\end{equation}
By (\ref{e46})-(\ref{e47}), we obtain
$$
\mathcal{A}(Jx)^{m-1}=\lambda (Jx)^{[m-1]},
$$
which means $Jx$ is also an H-eigenvector corresponding to $\lambda$. By assumptions, it follows that
$Jx =ax$ for some nonzero real constant, and $a$ is also an eigenvalue of
$J$. Then $a = ¡À1$. Therefore, $Jx = ¡Àx$, which implies that $x$ is either symmetric or skew-symmetric.
\qed

\setcounter{equation}{0}
\section{Centrosymmetric Cauchy tensors}

In \cite{chen14}, Chen and Qi introduced Cauchy tensors, and gave sufficient and necessary conditions for positive definiteness and semi-definiteness of even order Cauchy tensors.
In this section, we study centrosymmetric Cauchy tensors and give several sufficient and necessary conditions for a Cauchy tensor to be centrosymmetric. Furthermore, we prove that there are no odd dimension skew centrosymmetric Cauchy tensors. When the order is even, we prove that
all H-eigenvalues corresponding to nonzero H-eigenvalues of a centrosymmetric Cauchy tensor are symmetric vectors.    When the order is odd, we prove that the absolute
vectors of these H-eigenvalues are symmetric vectors.

Now, we first state the definition of Cauchy tensors.

\bd\label{def51}$^{\cite{chen14}}$ Let vector $c=(c_1,c_2,\cdots,c_n)\in \mathbb{R}^n$. Suppose that a real tensor $\mathcal{C}=(c_{i_1i_2\cdots i_m})$ is defined by
\begin{equation}\label{e51}
c_{i_1i_2\cdots i_m}=\frac{1}{c_{i_1}+c_{i_2}+\cdots+c_{i_m}},\quad j\in [m],~i_j \in [n].
\end{equation}
Then, we say that $\mathcal{C}$ is an order $m$ dimension $n$ symmetric Cauchy tensor and the vector $c\in \mathbb{R}^n$ is called the generating vector of $\mathcal{C}$.
\ed

In the sequence, a  centrosymmetric symmetric Cauchy tensor is called a centrosymmetric Cauchy tensor for simplicity.
\bt\label{thm51} Assume $\mathcal{C}$ is a Cauchy tensor defined as in (\ref{e51}). Let $c\in \mathbb{R}^n$ be the generating vector of $\mathcal{C}$.
Then Cauchy tensor $\mathcal{C}$ is centrosymmetric if and only if $c$ is symmetric i.e. $Jc=c$.
\et
\proof For sufficient conditions, suppose $c=(c_1,c_2,\cdots, c_n)$. By the definition of symmetric vectors, we have
$$c_i=c_{n-i+1},~i\in [n].$$
So, for any $i_1,i_2,\cdots,i_m \in [n]$, one has
$$
\begin{array}{rl} c_{i_1,i_2,\cdots,i_m}=&\frac{1}{c_{i_1}+c_{i_2}+\cdots+c_{i_m}} \\
=& \frac{1}{c_{n-i_1+1}+c_{n-i_2+1}+\cdots+c_{n-i_m+1}} \\
=& c_{n-i_1+1 n-i_2+1 \cdots n-i_m+1},
\end{array}
$$
which implies that $\mathcal{C}$ is a centrosymmetric Cauchy tensor.

For necessary conditions, assume Cauchy tensor $\mathcal{C}$ is centrosymmetric. Then we have
$$ c_{ii\cdots i}=c_{n-i+1 n-i+1\cdots n-i+1},~i\in [n],
$$
which means
$$\frac{1}{mc_i}=\frac{1}{mc_{n-i+1}},~i\in [n].$$
Thus $c_i=c_{n-i+1},~i\in [n]$ and $c$ is a symmetric vector.
\qed

\bt\label{thm52} Assume $\mathcal{C}$ is a Cauchy tensor defined as in (\ref{e51}). Then $\mathcal{C}$ is centrosymmetric if and only if
$$J\mathcal{C}=\mathcal{C}.$$
\et

\proof Let $c$ be the generating vector of $\mathcal{C}$. For any $i_1,i_2,\cdots,i_m \in [n]$, we have
\begin{equation}\label{e52}
\begin{array}{rl} (J\mathcal{C})_{i_1,i_2,\cdots,i_m}=&\sum_{t\in [n]}J_{i_1 t}c_{t i_2\cdots i_m} \\
=&c_{n-i_1+1 i_2\cdots i_m}\\
=&\frac{1}{c_{n-i_1+1}+c_{i_2}+\cdots+c_{i_m}}.
\end{array}
\end{equation}
Thus, from (\ref{e52}), we obtain that
$$J\mathcal{C}=\mathcal{C}$$
if and only if $c_i=c_{n-i+1},~i\in [n]$ i.e. $c$ is symmetric. By Theorem 5.1, we know that $J\mathcal{C}=\mathcal{C}$ if and only if
Cauchy tensor $\mathcal{C}$ is centrosymmetric.
\qed

By Theorem \ref{thm33}, we have the following result.
\bc\label{corol51} Assume $\mathcal{C}$ is a Cauchy tensor defined as in (\ref{e51}). Then $\mathcal{C}$ is centrosymmetric if and only if
$$\mathcal{C}J=\mathcal{C}.$$
\ec

\bt\label{thm53} Assume $\mathcal{C}$ is a Cauchy tensor defined as in (\ref{e51}).
Assume $n$ is even, then $\mathcal{C}$ is skew centrosymmetric if and only if $c$ is skew symmetric i.e. $Jc=-c$, where $c\in \mathbb{R}^n$
is the generating vector of $\mathcal{C}$.
\et
\proof When Cauchy tensor $\mathcal{C}$ is skew centrosymmetric, by Definitions \ref{def11} and \ref{def51}, we have
$$\frac{1}{mc_i}=c_{ii\cdots i}=-c_{n-i+1 n-i+1 \cdots n-i+1}=-\frac{1}{mc_{n-i+1}},~i\in [n].$$
Hence
$$c_i=-c_{n-i+1},~i\in [n],$$
which implies that $c$ is skew symmetric and the only if part holds.
For sufficient conditions, for any $i_1,i_2,\cdots,i_m \in [n]$, we have
$$\begin{array}{rl} c_{i_1i_2\cdots i_m}=&\frac{1}{c_{i_1}+c_{i_2}+\cdots+c_{i_m}}\\
=&-\frac{1}{c_{n-i_1+1}+c_{n-i_2+1}+\cdots+c_{n-i_m+1}}\\
=&-c_{n-i_1+1 n-i_2+1 \cdots n-i_m+1},
\end{array}
$$
where the second equality uses the fact that $c$ is skew symmetric. Thus Cauchy tensor $\mathcal{C}$ is skew centrosymmetric.
\qed

Here, it should be noted that there is no odd dimension skew centrosymmetric Cauchy tensor. If $\mathcal{C}$ is skew centrosymmetric
 Cauchy tensor and suppose $n$ is odd, let $i=\frac{n+1}{2}$, then
$$c_{ii\cdots i}=\frac{1}{mc_i}=-c_{n-i+1 n-i+1 \cdots n-i+1}=-\frac{1}{mc_{n-i+1}}=-\frac{1}{mc_i}.$$
Thus
$$\frac{1}{mc_i}=0,$$
which is a contradiction.

\bt\label{thm54} Assume order $m$ dimension $n$ Cauchy tensor $\mathcal{C}$ is defined as in (\ref{e51}). Let $c=(c_1,c_2,\cdots,c_n)$
be the generating vector of $\mathcal{C}$.
Suppose $\mathcal{C}$ is centrosymmetric. Then, for any H-eigenvector $x\in \mathcal{R}^n$ of $\mathcal{C}$ corresponding to a nonzero H-eigenvalue,
$x$ is symmetric when $m$ is even; $|x|$ is symmetric when $m$ is odd.
\et
\proof Since Cauchy tensor $\mathcal{C}$ is centrosymmetric, by Theorem \ref{thm51}, the generating vector $c$ is symmetric.
Suppose $x$ is an H-eigenvector of $\mathcal{C}$ corresponding to a nonzero H-eigenvalue $\lambda$. By Definition \ref{def22},
for $i\in [n]$, we have
$$
\begin{array}{rl} \lambda x_i^{m-1}=&(\mathcal{C}x^{m-1})_i\\
=&\sum_{i_2,i_3,\cdots,i_m \in [n]}c_{ii_2\cdots i_m}x_{i_2}x_{i_3}\cdots x_{i_m}\\
=&\sum_{i_2,i_3,\cdots,i_m \in [n]}\frac{1}{c_{i}+c_{i_2}+\cdots+c_{i_m}}x_{i_2}x_{i_3}\cdots x_{i_m}\\
=&\sum_{i_2,i_3,\cdots,i_m \in [n]}\frac{1}{c_{n-i+1}+c_{i_2}+\cdots+c_{i_m}}x_{i_2}x_{i_3}\cdots x_{i_m}\\
=&\lambda x^{m-1}_{n-i+1}.
\end{array}
$$
When $m$ is even, it holds that $x_i=x_{n-i+1}$, $i \in [n]$, which implies that $x$ is symmetric. When $m$ is odd, we
obtain $|x_i|=|x_{n-i+1}|$, $i \in [n]$, which implies that $|x|$ is symmetric.
\qed

\section{Final Remarks}

In this article, properties of centrosymmetric tensors and skew centrosymmetric tensors
are discussed. Some interesting results are natural extensions of the matrix case such as the products of centrosymmetric tensors, the
sufficient and necessary conditions for a tensor to be centrosymmetric and skew centrosymmetric. Spectral properties about H-eigenvalues and
H-eigenvectors of these tensors are also discussed. Furthermore, some symmetry properties of H-eigenvectors corresponding to nonzero H-eigenvalues of centrosymmetric Cauchy tensors are presented.  Some further questions are as follows.

{\bf Question 1}. How about the positive definiteness property of centrosymmetric tensors? Can we give some sufficient conditions just like the matrix case
in \cite {Andrew73}?

{\bf Question 2}. What are the properties of H-eigenvectors of skew centrosymmetric Cauchy tensors?


\begin{thebibliography}{99}
\bibitem{Andrew73} A L. Andrew, {\it Eigenvectors of certain matrices}, Linear Algebra Appl., 7 (1973) 151-162.
\bibitem{Anderson73} B.D.O. Anderson, E.I. Jury, {\it A simplified Schur-Cohn test}, IEEE Trans. Autom. Control, 18 (1973) 157-163.
\bibitem{Bu14} C. Bu, X. Zhang, J. Zhou, W. Wang, Y. Wei, {\it The inverse, rank and product of tensors}, Linear Algebra Appl. 446 (2014) 269-280.
\bibitem{Cant76} A. Cantoni, P. Butler, {\it Eigenvalues and eigenvectors of symmetric centrosymmetric matrices}
 Linear Algebra Appl., 13 (1976) 275-288.
\bibitem{chen14} H. Chen, L. Qi, {\it Positive Definiteness and Semi-Definiteness of Even Order Symmetric Cauchy Tensors},
arXiv preprint arXiv:1405.6363, 2014.
\bibitem{Chen14} Z. Chen, L. Qi, {\it Circulant Tensors with Applications to Spectral
Hypergraph Theory and Stochastic Process}, arXiv preprint
arXiv:1312.2752v7, 2014.
\bibitem{Demm97} J.W. Demmel, {\it Applied Numerical Linear Algebra}, SIAM, Philadelphia, 1997.
\bibitem{Ding13} W. Ding, L. Qi, Y. Wei, {\it M-Tensors and Nonsingular M-Tensors}, Linear Algebra Appl., 439 (2013) 3264-3278.
\bibitem{Ding14} W. Ding, L. Qi, Y. Wei, {\it Fast Hankel Tensor-Vector Products and Application to Exponential Data Fitting},
 arXiv preprint arXiv:1401.6238, 2014.
\bibitem{HH14} J. He, T.Z. Huang, {\it Inequalities for M-tensors},
J. Inequ. Appl.,  2014:114.
\bibitem{Muir83} M. Thomas, {\it A Treatise on the Theory of Deteminunts}, Dover, 1966 (originally published 1883).
\bibitem{Qi05} L. Qi, {\it Eigenvalue of a real supersymmetric tensor}, J. Symbol. Comput. 40 (2005) 1302-1324.
\bibitem{Qi14} L. Qi, {\it Hankel tensors: Associated Hankel matrices and Vandermonde decomposition}, Commun. Math. Sci., 12 (2014).
\bibitem{qi14} L. Qi, Y. Song, {\it An even order symmetric B tensor is positive definite},
Linear Algebra Appl., 457 (2014) 303-312.
\bibitem{QXX} L. Qi, C. Xu, Y. Xu, {\it Nonnegative tensor factorization, completely positive tensors and an hierarchical elimination algorithm}, to appear in: SIAM J. Matrix Anal.
Appl.
\bibitem{Shao13} J.Y. Shao, {\it A general product of tensors with applications},  Linear Algebra Appl., 439 (2013) 2350-2366.
\bibitem{Shao132} J.Y. Shao, H.Y. Shan, L. Zhang, {\it On some properties of the determinants of tensors}, Linear Algebra Appl., 439 (2013) 3057-3069.
\bibitem{Song14} Y. Song, L. Qi, {\it Some properties of infinite and finite dimension Hilbert tensors}, Linear Algebra Appl., 451 (2014) 1-14.
\bibitem{song14} Y. Song, L. Qi, {\it Properties of some classes of structured tensors}, to appear in: J. Optim. Theory Appl.
\bibitem{Weaver85} J.R. Weaver, {\it Centrosymmetric (cross-symmetric) matrices, their basic properties, eigenvalues, and eigenvectors}, American Mathematical Monthly, (1985)711-717.
\bibitem{YY14} P. Yuan, L. You, {\it Some remarks on P, P$_0$, B and B$_0$
tensors}, arXiv preprint arXiv:1405.1288, 2014.
\bibitem{Zhang12} L. Zhang, L. Qi, G. Zhou, {\it M-tensors and some applications},  SIAM J. Matrix Anal.
Appl., 35 (2014) 437-452.
\bibitem{Zhao14} X. Zhao, Q. Yang, {\it The spectral radius of nonnegative centrosymmetric tensor}, J. High School Numer. Math.,
36 (2014) 58-66.
\end{thebibliography}
\end{document}